\documentclass[preprint,3p]{elsarticle}
\usepackage{graphicx}
\usepackage{amsmath}
\usepackage{amssymb}
\usepackage{amsthm}
\usepackage{mathrsfs}
\usepackage{tikz-cd}
\usepackage{parskip}
\usepackage{hyperref}

\newtheorem{theorem}{Theorem}

\newtheorem{corollary}[theorem]{Corollary}

\newtheorem{lemma}[theorem]{Lemma}

\newtheorem{proposition}[theorem]{Proposition}

\theoremstyle{definition}
\newtheorem{algorithm}[theorem]{Algorithm}
\newtheorem{definition}[theorem]{Definition}
\newtheorem{example}[theorem]{Example}
\theoremstyle{remark}
\newtheorem{remark}[theorem]{Remark}

\newcommand{\closure}[2][3]{%
  {}\mkern#1mu\overline{\mkern-#1mu#2}}

\begin{document}

\title{An Algorithm For Checking Injectivity Of Specialization Maps From Elliptic Surfaces}
\author[1]{Tyler Raven Billingsley\corref{cor1}%
}
\address[1]{Department of Mathematics, Statistics and Computer Science, St. Olaf College, 
Northfield, MN 55057}
\ead{billin2@stolaf.edu}
\cortext[cor1]{Corresponding author}

\begin{abstract}
Let $E/\mathbb Q(t)$ be an elliptic curve and let $t_0 \in \mathbb Q$ be a rational number for which the specialization $E_{t_0}$ is an elliptic curve. Given a subgroup $M$ of $E(\mathbb Q(t))$ with mild conditions and $t_0 \in \mathbb Q$ coming from a relatively large subset $S_M \subset \mathbb Q$, we provide an algorithm that can show that the specialization map $\sigma_{t_0} : E(\mathbb Q(t)) \to E_{t_0}(\mathbb Q)$ is injective when restricted to $M$. The set $S_M$ is effectively computable in certain cases, and we carry out this computation for some explicit examples where $E$ is given by a Weierstrass equation.
\end{abstract}

\begin{keyword}
\MSC[2020]{11G05} \sep  Elliptic curves \sep Elliptic surfaces \sep Specialization \sep Sage
\end{keyword}

\maketitle

\section{Introduction}

Let $C$ be a (complete nonsingular) curve defined over a number field $k$ with function field $k(C)$. 
Let $E/k(C)$ be an elliptic curve defined by the Weierstrass equation $$y^2 = x^3 +A(t)x+B(t), \qquad A(t), B(t) \in k(C).$$ 
For any $t_0 \in C(k)$ such that the discriminant $4A(t)^3+27B(t)^2$ of $E$ does not vanish or have a pole at $t_0$, we define the elliptic curve $E_{0}/k$ using the Weierstrass equation $y^2 = x^3 + A(t_0)x+B(t_0)$. 
The specialization map at $t_0$ is the map 
$$\sigma_{t_0} : E(k(C)) \to E_{t_0}(k)$$
which takes the point $(x(t), y(t)) \in E(k(C))$ to $(x(t_0),y(t_0)) \in E_{t_0}(k)$. 
It is in fact a group homomorphism; that is, it respects the standard chord-and-tangent group laws on the domain and codomain. It is natural to ask what information can be extracted about the relationship between $E(k(C))$ and $E_{t_0}(k)$ through this homomorphism. 

In 1952, N\'eron proved a theorem regarding specialization of which the following result is a special case.

\begin{theorem} \cite{Neron-Spec} Let $k$ be a number field and let $E/k(t)$ be an elliptic curve. Then for infinitely many $t_0 \in k$ the specialization map
$$\sigma_{t_0} : E(k(t)) \to E_{t_0}(k)$$
is injective.
\end{theorem}

Thirty years later, Silverman \cite{Silverman-Spec} improved on N\'eron's result by proving that all but finitely many maps $\sigma_{t_0}$ are injective by showing that the set of $t_0 \in k$ for which injectivity fails is a set of bounded height.

While Silverman's result is effective, it might not give a practical way to find this finite set. So the aim of this article is to examine the injectivity of specialization homomorphisms in a different way. In particular, given an elliptic curve $E/k(t)$ given by a Weierstrass equation as above, how might one go about effectively determining the subset $\Sigma$ of $k$ containing those $t_0$'s for which the corresponding specialization maps are injective? We will focus on the case $k = \mathbb Q$ for concreteness. Gusi\'c and Tadi\'c \cite{GusicTadic} gave a method based on 2-descent which addresses this problem for curves with nontrivial $\mathbb Q(t)$-rational 2-torsion, but as they remark in their article the method does not immediately generalize to curves with trivial $\mathbb Q(t)$-rational 2-torsion. 

We consider a different method to find $t_0$'s for which the specialization map is injective. More specifically, we will discuss an approach to answering the following two questions.

	\begin{enumerate}
		\item Given some $t_0 \in \mathbb Q$, how can we effectively determine whether or not the specialization map at $t_0$ is injective?
		\item What is the set $\Sigma$ of $t_0 \in \mathbb Q$ such that the specialization map at $t_0$ is injective?
	\end{enumerate}

We will approach these questions from the perspective of irreducibility. We review the classical N\'eron specialization theorem to relate the above questions to the question of irreducibility of polynomials after specialization, and we give an algorithm that can be used to find a Hilbert set which intersects $\Sigma$ in an infinite set $S$. The approach taken here was inspired by a mathoverflow post by Silverman \cite{Silverman-MOpost}. The algorithm is effective in certain cases, and we carry out the algorithm on some examples. In doing so, the elliptic curve packages in Sage \cite{sagemath} and Magma \cite{MR1484478} were indispensible. They will be named wherever they were used.

\section{Preliminaries} \label{preliminaries}
We begin by reviewing the basics on division polynomials and Hilbert sets that are used in the proof of N\'eron's specialization theorem. 

Let $K$ be a field of characteristic zero. Given a $K$-rational point on some elliptic curve $E$, finding a point $Q$ with $nQ = P$ requires solving polynomial equations over $K$. The polynomials which appear in this way are called division polynomials. Fix $A, B \in K$ with $4A^3 + 27B^2 \neq 0$ and define
$$\mathcal O = K[x,y]/(y^2-(x^3+Ax+B));$$
that is, $\mathcal O$ is the coordinate ring of the elliptic curve
$$E : y^2 = x^3+Ax+B$$
defined over $K$. Define the following sequence of polynomials in $\mathcal O$.
\begin{align*}
\psi_0 &= 0, \\
\psi_1 &= 1, \\
\psi_2 &= 2y, \\
\psi_3 &= 3x^4+6Ax^2+12Bx - A^2, \\
\psi_4 &= 4y(x^6+5Ax^4 + 20Bx^3 -5A^2x^2-4ABx-8B^2-A^3), \\
&\vdots \\
\psi_{2m+1} &= \psi_{m+2}\psi_m^3-\psi_{m-1}\psi^3_{m+1} \text{ for } m \geq 2, \\
\psi_{2m} &= \left(\frac{\psi_m}{2y}\right) \cdot (\psi_{m+2}\psi^2_{m-1}-\psi_{m-2}\psi^2_{m+1}) \text{ for } m \geq 3.
\end{align*}
We have the following standard facts about the polynomials above, which we state without proof.
\begin{itemize}
	\item The polynomials $\psi_{2n+1}$, $\psi_{2n}/y$ and $\psi_{2n}^2$ depend only on $x$.
	\item For $n \geq 1$, the set of roots of $\psi_{2n+1}$ is the set of $x$-coordinates of the nonzero $(2n+1)$-torsion points.
	\item For $n \geq 2$, the set of roots of $\psi_{2n}/y$ is the set of $x$-coordinates of the nonzero $(2n)$-torsion points which are not 2-torsion. Additionally, the set of roots of $x^3+Ax+B$ is the set of $x$-coordinates of nonzero 2-torsion points. Since $y^2 = x^3+Ax+B$, we see that the set of roots of $\psi_{2n}^2$ is the set of $x$-coordinates of all nonzero $(2n)$-torsion points. However, it is worth noting that this polynomial is not separable.
	\item Combining the previous two statements, for $n \geq 2$ we have that the set of roots of $\psi_n^2$ is the set of $x$-coordinates of the nonzero $n$-torsion points.
	\item If we additionally define
$$\phi_n = x\psi_n^2 - \psi_{n+1}\psi_{n-1}, \qquad \omega_n = \frac{\psi_{n+2}\psi^2_{n-1}-\psi_{n-2}\psi^2_{n+1}}{4y},$$
then for any $Q = (x_Q, y_Q) \in E(K) \setminus E[n](K)$ we have that 
\begin{equation} \label{eqn:mult_div_poly} nQ = \left(\frac{\phi_n(x_Q)}{\psi^2_n(x_Q)},\frac{\omega_n(x_Q,y_Q)}{\psi_n^3(x_Q,y_Q)}\right).
\end{equation}
Note that, if $Q \in E[n](K)$, by the above discussion we must have $\psi^2_n(x_Q) = 0$.
	\item The polynomials $\psi_n^2$ and $\phi_n$ have degrees $n^2-1$ and $n^2$, respectively.
\end{itemize}

Using Equation \eqref{eqn:mult_div_poly}, the point $P \neq O$ is divisible by $n$ in $E(K)$ if and only if there exists a point $Q = (x_Q,y_Q) \in E(K)$ such that
$$x_P = \frac{\phi_n(x_Q)}{\psi_n^2(x_Q)}, \qquad y_P = \frac{\omega_n(x_Q,y_Q)}{\psi_n^3(x_Q,y_Q)}.$$
Focusing on the equation for $x_P$, we must have that
$$\phi_n(x_Q) - x_P\psi^2_n(x_Q) = 0.$$
Thus $x$-coordinates of points $Q$ with $nQ = P$ satisfy the polynomial
\begin{equation} \label{eqn:n_div_poly} d_{n,P}(x) =  \phi_n(x) - x_P\psi^2_n(x) = 0. \end{equation}
We call $d_{n,P}(x)$ the \textbf{$n$-division polynomial of the point $P$}.
\begin{lemma} \label{lem:div_poly_of_P}
Let $d_{n,P}(x)$ be the $n$-division polynomial of a point $P \in E(K) \setminus E[2](K)$. Then $d_{n,P}(x)$ has a root in $K$ if and only if $P$ is divisible by $n$ in $E(K)$. 
\end{lemma}
\begin{proof}
If $P$ is divisible by $n$ in $E(K)$, say $nQ = P$, use Equation \eqref{eqn:mult_div_poly} and repeat the above derivation to show that $d_{n,P}(x)$ has a root in $K$ which is the $x$-coordinate of the point $Q$. Conversely, suppose $d_{n,P}(x_Q) = 0$ for some $x_Q \in K$. Since $d_{n,P}(x)$ is a polynomial of degree $n^2$ and there are $n^2$ $n$-division points of $P$ in $\bar K$ all with distinct $x$-coordinates (since $P$ is not 2-torsion), we have some $y_Q \in \bar K$ such that the point $Q = (x_Q,y_Q) \in E(\bar K)$ has the property that $nQ = P$. Now, for any $\tau \in \text{Gal} (\bar K/K)$, we have that
$$n \cdot \tau(Q) = P$$
with $x_{\tau(Q)} = x_Q$. Using the Weierstrass equation we then have $y_{\tau(Q)} = \pm y_Q$, so $\tau(Q) = \pm Q$. Since $P$ is not 2-torsion we have that $P \neq -P$. Therefore we cannot have $\tau(Q) = -Q$, because otherwise $n\cdot \tau(Q) = n(-Q) = -nQ = -P \neq P$. Hence $\tau(Q) = Q$ for every $\tau \in \text{Gal} (\bar K/K)$. Thus $y_Q \in K$. 

\begin{remark} \label{rem:2_torsion_bad}
It is possible that $P \in E[2](K)$ and $d_{n,P}(x)$ has $K$-rational roots, but $P$ is not divisible by $n$ in $E(K)$. For example, let $n=2$ and $K = \mathbb Q$. Set
$$E: y^2 = x^3 + 503844x - 45019744.$$
Then $P = (88,0) \in E[2](\mathbb Q)$ with
$$d_{2,P}(x) = ((x - 814)(x + 638))^2,$$
but neither $814$ nor $-638$ are $x$-coordinates of points in $E(\mathbb Q)$ (but they are $x$-coordinates of points $Q \in E(\bar{\mathbb Q})$ with $2Q = P$).
\end{remark}
\end{proof}

We now mention some basic facts regarding Hilbert sets. For more details about Hilbert sets and Hilbert's irreducibility theorem in greater generality, see Lang \cite{Lang-DiophantineGeometry}. 

\begin{definition}
Let $f_1(t,x), ... ,f_n(t,x) \in \mathbb Q[t,x]$ be irreducible polynomials over $\mathbb Q$. The {\upshape Hilbert subset of} $\mathbb Q$ corresponding to the $f_i$'s is the set of all $t_0 \in \mathbb Q$ such that each $f_i(t_0,x) \in \mathbb Q[x]$ is irreducible over $\mathbb Q$. That is, the Hilbert subset is the set of $t_0$'s for which all the polynomials remain irreducible upon specialization at $t = t_0$.
\end{definition}

Notice that the intersection of any two Hilbert sets is a Hilbert set - it corresponds to the union of the sets of polynomials defining the two Hilbert sets. Additionally, note that the preceding definition varies slightly from that in \cite{Lang-DiophantineGeometry}, where a Hilbert set is also allowed to intersect finitely many Zariski open sets of $\mathbb A^1_{\bar{\mathbb Q}}$ - that is, they might differ from a Hilbert set as above by finitely many points. Although we will have a use for such a set later, we prefer to use the term ``Hilbert set" for sets precisely corresponding to irreducible specializations.

\begin{theorem}[Hilbert's Irreducibility Theorem]
Every Hilbert subset of $\mathbb Q$ is infinite.
\end{theorem}

\begin{proof}
\cite[Chapter 9 \S 2 Corollary 2.5]{Lang-DiophantineGeometry}.
\end{proof}

Recall that the natural density of a subset $T$ of $\mathbb N$ is the limit
$$\lim_{n \to \infty} \frac{\#\{k \in T \mid k \leq n\}}{n}.$$

\begin{proposition}  \label{prop:density1}
For any Hilbert subset $H$ of $\mathbb Q$, $H \cap \mathbb N$ has natural density 1.
\end{proposition}

\begin{proof}
Lang \cite[Chapter 9 \S 2 Corollary 2.3]{Lang-DiophantineGeometry} states that, for $n$ large enough, 
$$n - n^\alpha \leq \#\{k \in H \cap \mathbb N \mid k \leq n \} $$
for some fixed $\alpha$ with $0 < \alpha < 1$ independent of $n$. Hence for $n$ large enough, we have
$$1-\frac{1}{n^{1-\alpha}} \leq \frac{\#\{k \in H \cap \mathbb N \mid k \leq n \}}{n} \leq 1,$$
so 
$$\lim_{n \to \infty} \frac{\#\{k \in H \cap \mathbb N \mid k \leq n\}}{n} = 1.$$
\end{proof}

\section{N\'eron's Specialization Theorem} \label{sec:neron}

N\'eron's specialization theorem (for elliptic curves) is the result of applying Hilbert's irreducibility theorem to division polynomials. Following Serre \cite[Chapter 11]{Serre-MordellWeil}, we start by recalling a completely group-theoretic fact.

\begin{proposition} \label{prop:grouptheory}
Let $n$ be a positive integer and let $\phi: M \to N$ be a homomorphism of abelian groups with the following properties.
	\begin{enumerate}
	\item $M$ is finitely generated.
	\item The induced map $\bar\phi: M/nM \to N/nN$ is injective.
	\item $\phi|_{M[n]}$ gives an isomorphism $M[n] \cong N[n]$.
	\item $\phi|_{M_\text{tors}}$ is injective.
	\end{enumerate}
Then $\phi$ is injective.
\end{proposition}

Fix an elliptic curve $E$ over $\mathbb Q(t)$, and set $\phi = \sigma_{t_0}$ the specialization homomorphism for a fixed $t_0 \in \mathbb Q$, $M = E(\mathbb Q(t))$, $N = E_{t_0}(\mathbb Q)$ and a positive integer $n \geq 2$. Then conditions 1 and 4 above are always true; indeed, condition 1 is the function field version of the Mordell-Weil Theorem \cite[Chapter III \S 6]{Silverman2}, and condition 4 follows from basic results on formal groups of elliptic curves and their relationship to reduction mod $p$ found in Silverman \cite[Chapter VII \S 3]{Silverman1}. It is true that conditions 2 and 3 hold inside of a Hilbert set, but instead of proving this fact directly we will replace $E(\mathbb Q(t))$ with specific subgroups (which will include $E(\mathbb Q(t))$ itself) and prove a more general statement.

First, notice that condition 2 is equivalent to the following statement.
\begin{equation} \label{condition2} \text{For any } a \in M \text{ such that } \phi(a) \text{ is divisible by } n \text{ in } N, a \text{ is divisible by } n \text{ in } M.
\end{equation}
Now suppose $M < E(\mathbb Q(t))$ and $\phi = \sigma_{t_0}|_M$. If we have some $a \in M$ with $\phi(a)$ divisible by $n$, then even if we are able to conclude that $a$ is divisible by $n$ in $E(\mathbb Q(t))$, say $a = nb$ for some $b \in E(\mathbb Q(t))$, in order to verify \eqref{condition2} we would be required to check that $b \in M$. Unfortunately, this is a difficult problem. To avoid this, we will restrict ourselves to subgroups $M$ in which whenever $a \in M$ has the property that $a = nb$ for some $b \in E(\mathbb Q(t))$, we always know that $b \in M$. That is, condition 2 of Proposition \ref{prop:grouptheory} is satisfied for the inclusion $M \to E(\mathbb Q(t))$. Of course, this isn't always true; for example, if a subgroup $G < E(\mathbb Q(t))$ contains an element of infinite order then $M = nG$ always fails to have this property (recall that $G$ is finitely generated). Regarding condition 3, since specialization is always injective on torsion, note that if $M[n]$ is a proper subgroup of $E[n](\mathbb Q(t))$, then $\phi(M[n])$ is also always a proper subgroup of $E_{t_0}[n](\mathbb Q)$. Hence we also make the minor additional assumption that $M[n] = E[n](\mathbb Q(t))$. Thus we will assume all hypotheses of Proposition \ref{prop:grouptheory} for the inclusion \break $M \to E(\mathbb Q(t))$. 

Now suppose that we have $M < E(\mathbb Q(t))$ with the inclusion satisfying the hypotheses of Proposition \ref{prop:grouptheory}, $N = E_{t_0}(\mathbb Q)$ and $\phi = \sigma_{t_0}$. In our goal of finding a Hilbert set on which conditions 2 and 3 of Proposition \ref{prop:grouptheory} hold for these very specific $M, N, \phi$, one might hope that we could discard either condition 2 or 3 and still maintain the conclusion of Proposition \ref{prop:grouptheory}. Unfortunately neither condition implies the other, as the following two examples illustrate. We first look at an example where condition 2 holds but condition 3 does not.

\begin{example} \label{ex:torsion_needed}
Let $E: y^2 = x^3-(t^2+27)x+(10t^2+48t+90)$, $\phi = \sigma_{30}$, $M = E(\mathbb Q(t))$, $N = E_{30}(\mathbb Q)$ and $n = 2$. In \cite{Shioda-Examples}, this elliptic curve is shown to have Mordell-Weil rank 4 and to have no nontrivial torsion points over $\mathbb{Q}(t)$. The Mordell-Weil group is generated by the four points $P_1 = (9, t+24), P_2 = (6,2t+12), P_3 = (1,3t+8)$ and $P_4 = (t+3,4t+6)$, so a complete set of representatives for the nonidentity cosets of $2E(\mathbb Q(t))$ in $E(\mathbb Q(t))$ is 
$$\left\{\sum_{i \in C} P_i \mid C \subset \{1,2,3,4\}, C \neq \emptyset\right\}.$$
One can check (for instance, using the {\upshape \texttt{EllipticCurve}} method {\upshape \texttt{division\_points(2)}} in Sage) that for $t=30$ the specialization of each of these 15 points is not divisible by 2 in $E_{30}$. Thus $\closure[1]{\sigma_{30}}: E(\mathbb Q(t))/2E(\mathbb Q(t)) \to E_{30}(\mathbb Q)/2E_{30}(\mathbb Q)$ is injective. However, the Mordell-Weil group of $E_{30}$ is $\mathbb{Z}^3 \times \mathbb{Z}/2\mathbb{Z}$, so $\sigma_{30}$ cannot be injective. In particular, condition 3 of Proposition \ref{prop:grouptheory} does not hold, but condition 2 does hold.
\end{example}

Next, we consider an example where condition 3 holds but condition 2 does not.

\begin{example}
Let $E: y^2 = x^3-t^2x+t^2, \phi = \sigma_{2}$, $M = E(\mathbb Q(t))$, $N = E_{2}(\mathbb Q)$ and $n = 2$. One can check that $E(\mathbb Q(t)) \cong \mathbb Z^2$ with generators
$$P = (t,t), Q = (0,t);$$
see \S \ref{sec:examples}. Then, using Sage, one can check that the specialization $E_2$ has $E_2(\mathbb Q) \cong \mathbb Z$. Hence $\closure[1]{\sigma_2}$ is a map from a group of order 4 to a group of order 2, so $\closure[1]{\sigma_2}$ cannot be injective.
\end{example}

Thus we will need to check both of the conditions as part of the following proof of N\'eron's specialization theorem for subgroups. For a more general version of this theorem, see \cite{Lang-DiophantineGeometry}.

\begin{theorem} \label{thm:neron-spec}
Let $E/\mathbb Q(t)$ be a nonconstant elliptic curve given by the Weierstrass equation
$$y^2 = x^3+A(t)x+B(t)$$
and let $M < E(\mathbb Q(t))$ be a subgroup of rank at least 1 such that the inclusion $M \to E(\mathbb Q(t))$ satisfies the hypotheses of Proposition \ref{prop:grouptheory}. Then there exists a set $S_M$ which differs from a Hilbert set by finitely many elements such that for each $t_0 \in S_M$ the specialization map $\sigma_{t_0}|_M: M \to E_{t_0}(\mathbb Q)$ is injective.
\end{theorem}

\begin{proof}

By the preceding comments, it suffices to show that the set of $t_0 \in \mathbb Q$ for which $\phi = \sigma_{t_0}|_M$ satisfies conditions 2 and 3 of Proposition \ref{prop:grouptheory} differs from a Hilbert set by finitely many elements. Let $n \geq 2$ and suppose $P_1, ... , P_k$ are a set of representatives for the nonzero elements of $M/nM$. Since $M$ has rank at least 1, by possibly changing some of the $P_i$'s by an element of $nM$ we may assume that no $P_i$ is 2-torsion. Let $d_{n,P_i}(t,x)$ be the $n$-division polynomial of $P_i$ (\S \ref{preliminaries}). After clearing denominators, we may assume $d_{n,P_i}(t,x) \in \mathbb Q[t][x]$. Because condition 2 of Proposition \ref{prop:grouptheory} holds for the inclusion $M \to E(\mathbb Q(t))$, no $P_i$ is divisible by $n$ in $E(\mathbb Q(t))$ and thus, by Lemma \ref{lem:div_poly_of_P} each $d_{n,P_i}(t,x)$ has no roots in $\mathbb Q(t)$ as a polynomial in $x$ - that is, the irreducible factorization of $d_{n,P_i}(t,x)$ in $\mathbb Q[t][x]$ has no factors with ($x$-)degree 1. Let $H_1$ be the Hilbert set corresponding to all irreducible factors of all the $d_{n,P_i}$'s, then remove any rational number from $H_1$ which appears as a zero of a coefficient in an irreducible factor; call this set $S_1$. So for any $t_0 \in S_1$ each irreducible factor of each $d_{n,P_i}(t,x)$ remains irreducible upon specialization, and since none of the coefficients vanish the $x$-degree is preserved. Thus $d_{n,P_i}(t_0,x)$ has no roots in $\mathbb Q$. Since the roots of this polynomial are $x$-coordinates of points $Q_{t_0}$ such that $nQ_{t_0} = P_{i,t_0}$, $P_{i,t_0}$ is not divisible by $n$ in $E_{t_0}(\mathbb Q)$, and thus condition 2 is satisfied.

Next, using notation from \S \ref{preliminaries}, consider the polynomial $\psi_n^2$. Recall that this polynomial has the set of $x$-coordinates of the $n$-torsion points of $E(\closure[1]{\mathbb Q(t)})$ as its roots. Clearing denominators, we assume that $\psi_n^2 \in \mathbb Q[t][x]$. Let $\{r_1, ... , r_l\}$ be the $\mathbb Q(t)$-rational roots of $\psi_n^2$ which do not correspond to $\mathbb Q(t)$-rational $n$-torsion points (by Remark \ref{rem:2_torsion_bad} this is a possibility when $n$ is even), and let $f_1, ... ,f_l$ be the polynomials obtained by clearing the denominators in the expressions
$$x^2-(r_i^3+A(t)r_i+B(t)).$$
Notice that, since each $r_i$ is not the $x$-coordinate of a point in $E(\mathbb Q(t))$, we have that each $f_i$ is irreducible over $E(\mathbb Q(t))$. Let $H_2$ be the Hilbert set corresponding to the irreducible factors of $\psi_n^2$ of degree at least 2 and the polynomials $f_i$, and then remove any rational number from $H_2$ which appears as a zero of a coefficient; call this set $S_2$. Then, upon specialization at $t_0 \in S_2,$ each irreducible factor of degree at least 2 remains irreducible of degree at least 2, and the fact that the polynomials $f_i$ remain irreducible of $x$-degree 2 means that $r_i(t_0)^3+A(t_0)r_i(t_0)+B(t_0)$ is not a square in $\mathbb Q$ so that $r_i(t_0)$ is not the $x$-coordinate of a point in $E_{t_0}(\mathbb Q)$. Because of this, $E_{t_0}$ gains no new $\mathbb Q$-rational $n$-torsion points, so condition 3 is satisfied. Recall that specialization is injective on torsion.

Finally, remove from $S_1 \cap S_2$ the poles of $A$, the poles of $B$ and all $t_0$ such that $E_{t_0}$ is not smooth; call this set $S_M$. Then $S_M$ has the required property.

\end{proof}

\begin{remark}
The choice of $n \geq 2$ in the proof of Theorem \ref{thm:neron-spec} has no impact on the proof - any value of $n$ works. However, as $n$ becomes larger so does $M/nM$, so choosing $n$ to be smaller can help when working with examples.
\end{remark}

\begin{corollary} \label{cor:density1-spec}
Let $E$ and $M$ be as in Theorem \ref{thm:neron-spec}. Then the set $\Sigma_M$ of all $k_0 \in \mathbb N$ such that the specialization map $\sigma_{k_0}|_M$ is injective has density 1.
\end{corollary}
\begin{proof}
Proposition \ref{prop:density1}.
\end{proof}

We conclude the section with a summary of how the previous proof yields an algorithm that can often be used to check when a specialization map is injective. 

\begin{algorithm} \label{algorithm}
Let $E/\mathbb Q(t)$ be an elliptic curve given by a Weierstrass equation and let $M$ be a subgroup as in Theorem \ref{thm:neron-spec}. Fix $n \geq 2$, though in the examples to follow we always choose $n = 2$.
	\begin{enumerate}
	\item Let $\{ P_1, ..., P_k\}$ be a set of representatives of the nonzero cosets of $nM$ in $M$, taking care not to choose a $2$-torsion point.
	\item For each $P_i$, compute $d_{n,P_i}(t,x)$ and clear denominators to assume that $$d_{n,P_i}(t,x) \in \mathbb{Q}[t][x].$$
	\item Compute the collection of polynomials $f_i$ as in Theorem \ref{thm:neron-spec} and the non-linear irreducible factors of the division polynomial $\psi_n^2.$
	\item Compute the Hilbert set corresponding to the irreducible factors of the polynomials above, then compute the set $S_M$ as in Theorem \ref{thm:neron-spec} by removing the poles of $A$, the poles of $B$ and those $t_0$'s for which coefficients of at least one of the above polynomials vanish or $E_{t_0}$ is not smooth.
	\end{enumerate}
\end{algorithm}

\begin{remark} \label{rem:roots}
In practice, one only needs that the specialized polynomials have no roots, which is weaker than asking that all irreducible factors remain irreducible.
\end{remark}

\begin{remark} \label{rem:explicit}
Notice that step 4 requires checking that various irreducible factors remain irreducible upon specialization; for a fixed $t_0 \in \mathbb Q$, this can often be done by inspection or with computer software such as Sage.
\end{remark}

\begin{remark} \label{rem:generalized_algorithm} Let $E$ be an elliptic curve in Weierstrass form defined over the function field $K$ of a curve defined over a number field $k$, and suppose the curve is given by an explicit equation. While Algorithm \ref{algorithm} was written specifically for elliptic curves over $\mathbb Q(t)$, the above algorithm can be adjusted to work for specializing at $k$-rational points of the curve. In particular, it can be used to check injectivity of a specific specialization map, as in Remark \ref{rem:explicit}.

\end{remark}

\section{Examples Using the Irreducibility Algorithm} \label{sec:examples}

In this section, we discuss some explicit examples of utilizing Algorithm \ref{algorithm} with the modification mentioned in Remark \ref{rem:roots}. We consider two examples, one with a full Mordell-Weil group of rank 2 and another with a subgroup of rank 2. 

We first consider
$$E: y^2 = x^3 - t^2x + t^2.$$ 
Our goal is to find an infinite set of rational numbers for which the corresponding specialization maps (on all of $E(\mathbb Q(t))$) are injective. First, we need generators of $E(\mathbb Q(t))$ in order to use Algorithm \ref{algorithm}. Set $P = (t,t)$ and $Q = (0,t)$. One can check (using Magma) that the determinant of the canonical height matrix of $P$ and $Q$ is nonzero, and thus $P$ and $Q$ are linearly independent in $E(\mathbb Q(t))$. In addition, using Magma's implementation of Tate's algorithm \cite[Chapter IV \S 9]{Silverman2} and combining the resulting information with the Shioda-Tate formula \cite{Shioda-MWlattices}, the Mordell-Weil rank of $E/\closure[1]{\mathbb Q}(t)$ is 2. Hence the rank of $E(\mathbb Q(t))$ is 2. In order to show that $P$ and $Q$ generate $E(\mathbb Q(t))$, we will use specialization in a way that is motivated by (but different from) the method outlined in \cite{Stoll-Spec}. First, we show that $E(\mathbb Q(t))$ has trivial torsion. Since specialization is injective on torsion, it suffices to show that a single specialization has trivial torsion. To see this, consider the following example, which (at the same time) highlights Remark \ref{rem:explicit} and shows how injectivity of the specialization map for an individual $t_0 \in \mathbb Q$ can often be checked directly using a computer algebra system such as Sage.

\begin{example} \label{ex:check_2div_points}
Let $t_0 = 5$ and let $M = \langle P, Q \rangle$. Consider the specialized curve 
$$E_{5}: y^2 = x^3 - 25x + 25.$$  Note that
\begin{align*}
M/2M &= \{ 0, P, Q, P+Q \} \\
&= \{ 0, (t,t), (0, t), (-t, -t) \}.
\end{align*}
Run the following code in a Sage worksheet.
\begin{verbatim}
t = 5
Espec = EllipticCurve([-t^2,t^2])
Pspec = Espec(t,t)
Qspec = Espec(0,t)
print("The 2-division points of Pspec are: "
	    + str(Pspec.division_points(2)))
print("The 2-division points of Qspec are: "
	    + str(Qspec.division_points(2)))
print("The 2-division points of Pspec+Qspec are: "
	    + str((Pspec+Qspec).division_points(2)))
print("The torsion points of Espec are: "
	    + str(Espec.torsion_points()))
\end{verbatim}
The output is the following.
\begin{verbatim}
The 2-division points of Pspec are: []
The 2-division points of Qspec are: []
The 2-division points of Pspec+Qspec are: []
The torsion points of Espec are: [(0 : 1 : 0)]
\end{verbatim}
We interpret the output as follows. In the context of Proposition \ref{prop:grouptheory}, set $\phi = \sigma_5|_M$, $N = E_5(\mathbb Q)$, $n = 2$ and use $M$ as already defined. First, since $E_5(\mathbb Q)$ has no torsion and specialization is injective on torsion, we see that $E(\mathbb Q(t))$ (and thus $M$) also has no torsion. Thus condition 3 holds. Condition 2 is equivalent to the generators of $M/2M$ not being divisible by 2 in $E(\mathbb Q)$ upon specialization, which is shown by the above output. Hence condition 2 holds. Finally, we need to show that conditions 2 and 3 hold for the inclusion $M \to E(\mathbb Q(t))$. Clearly condition 3 holds since we've shown that $E(\mathbb Q(t))$ has no torsion, and we can show condition 2 by using Sage to show that $P, Q$ and $P+Q$ have no 2-division points in $E(\mathbb Q(t))$ using similar commands to those above. Thus the specialization map $\sigma_5|_M$ is injective.
\end{example}

It remains to show that $P$ and $Q$ generate $E(\mathbb Q(t))$. Consider the specialization at $t_0 = 5$ as in Example \ref{ex:check_2div_points}. Sage yields that $E_5$ has Mordell-Weil group $\mathbb Z^2$ over $\mathbb Q$ with generators $(-1,7)$ and $(0,5)$. Since $(-1,7) + (0,5) = (5,5),$ we may instead use $(5,5) = P_5$ and $(0,5) = Q_5$ as generators. Fix the bases $\{P,Q\}, \{P_5,Q_5\}$ for $M, E_5(\mathbb Q)$, respectively. After fixing some basis for $E(\mathbb Q(t))$ (which has 2 elements), we let the matrix $A$ represent the inclusion $M \to E(\mathbb Q(t))$ and the matrix $B$ represent the specialization map $\sigma_5$. We then have a sequence
$$M \overset{A}\to E(\mathbb Q(t)) \overset{B}\rightarrow E_5(\mathbb Q).$$
The composition $BA$ is the specialization map $\sigma_5|_M$. Since this maps generators of $M$ to generators of $E_5(\mathbb Q)$, $BA$ is the identity matrix. Hence $A$ is invertible, so the inclusion $M \to E(\mathbb Q(t))$ is surjective. Hence $E(\mathbb Q(t)) \cong \mathbb Z^2$ with generators $P = (t,t)$ and $Q = (0,t)$.

Before moving forward with using Algorithm \ref{algorithm} to find injective specialization maps, it is important to notice that success of this method for a fixed $n$ is not equivalent to injectivity of the specialization map. 
We can't hope for this to be true since it succeeds on (most of) a Hilbert set and Hilbert sets often have infinite complements, whereas Silverman's specialization theorem states that the specialization map fails to be injective for only finitely many rational numbers.
The next example illustrates the failure of this equivalence.

\begin{example} \label{ex:linear_algebra}
Let $t_0 = 27$. On the elliptic curve $E_{27}: y^2 = x^3 - 729x + 729$, notice that
$$[2](-9,81) = (27,27) = P_{27},$$
so our criterion (for $n=2$) cannot conclude that $\sigma_{27}$ is injective because condition 2 of Proposition \ref{prop:grouptheory} fails. A check using Sage shows that $E_{27}(\mathbb Q) \cong \mathbb Z^2$ with generators $R_1 = (-9,81)$ and $R_2 = (-27,27)$. Now $P_{27} = 2R_1$ and $Q_{27} = -(2R_1+R_2)$, meaning the matrix of the specialization map $\sigma_{27}$ with respect to the ordered bases $\{P,Q\}$ and $\{R_1,R_2\}$ is
$$\begin{bmatrix}
2 & -2 \\
0 & -1
\end{bmatrix}.$$
The determinant of this matrix is $-2 \neq 0$, so $\sigma_{27}$ is injective.
\end{example}

We now carry out Algorithm \ref{algorithm} for $n=2$. As in \S \ref{preliminaries}, we find the the polynomials in steps 3 and 4 to be
\begin{align*}
d_{2,P}(t,x) &= x^4+2t^2x^2-8t^2x+t^4-t(4x^3-4t^2x+4t^2), \\
d_{2,Q}(t,x) &= x^4+2t^2x^2-8t^2x+t^4, \\
d_{2,P+Q}(t,x) &= x^4+2t^2x^2-8t^2x+t^4+t(4x^3-4t^2x+4t^2), \text{ and} \\
g(t,x) &= x^3-t^2x+t^2,
\end{align*}
where $g(t,x) = \psi_2^2/4$. Notice that all four polynomials are irreducible over $\mathbb Q[t,x]$. As in Remark \ref{rem:roots}, we need to find $t_0$'s for which the specialized polynomials have no roots in $\mathbb Q$. Equivalently, we need to find $t_0$'s for which the curves defined by the polynomials have no rational points of the form $(t_0,x_0)$. Set
\begin{align*}
C_P: & \,\, d_{2,P}(t,x) = 0, \\
C_Q:&\,\, d_{2,Q}(t,x) = 0, \\
C_{P+Q}:& \,\, d_{2,P+Q}(t,x) = 0 \\
C_2:& \,\, g(t,x) = 0.
\end{align*}
Using Sage, all of the curves are rational over $\mathbb Q$ and have rational points, hence they have infinitely many rational points. Because of this, we will restrict to $t_0 \in \mathbb N$ and examine the case of specializing at natural numbers, which similarly reduces to looking at integral points on the curves (since each polynomial is monic in $x$). We first prove, using elementary methods, that the only obstruction to success of the method for $t_0 > 2$ comes from $C_{P}$. We begin with an algebraic lemma which will make analyzing $C_{P+Q}$ easy.

\begin{lemma} \cite{Rees-Quartic} \label{lemma:quartic}
Consider a (depressed) quartic polynomial 
$$p(x) = x^4+qx^2+rx+s \in \mathbb{Q}[x]$$
with discriminant $\Delta > 0$. If $q < 0$ and $s < q^2/4,$ then $p$ has four distinct roots in $\mathbb{C} \setminus \mathbb{R}$.
\end{lemma}

\begin{proposition} 
The curves $C_{P+Q}, C_Q$ and $C_2$ each have no integral points with $t_0 > 2$.
\end{proposition}

\begin{proof}
$C_{P+Q}$: Notice that the discriminant of $d_{2,P+Q}$ as a polynomial in $x$ is 
$$16384t^{10} - 110592t^8.$$
This is positive for $t_0 > 2$. The corresponding depressed quartic (in $x$) is 
$$x^4 - 4t^2x^2 - 8t^2x + 4t^4 + 12t^3.$$
 By Lemma $\ref{lemma:quartic}$, this quartic has no real roots in $x$ for $t_0 > 2$.

$C_Q$: Fix $t_0 \in \mathbb N$ and suppose $(t_0,x_0)$ is an integral point on $C_Q$. We make 3 cases based on possible values of $x_0$.
\begin{itemize}
\item[] Case 1: If $x_0 \leq 0$, each nonzero term of $d_{2,Q}(t_0,x_0)$ is positive. Thus $(t_0,x_0)$ is not a point on $C_Q$.
\item[] Case 2: Let $x_0 \geq 4$. Note that
$$d_{2,Q}(t_0,x_0) = x_0^4 + 2t_0^2x_0^2-8t_0^2x_0+t_0^4 \geq 256+32t_0^2-32t_0^2+t_0^4 = 256+t_0^4 > 0,$$
so $(t_0,x_0)$ is not a point on $C_Q$.
\item[] Case 3: Suppose $1 \leq x_0 \leq 3$. We have the following three polynomials in $t$:
$$d_{2,Q}(t,1) = t^4 -6t^2+1,$$
$$d_{2,Q}(t,2) = t^4-8t^2+16, \text{ and}$$
$$d_{2,Q}(t,3) = t^4-6t^2+81.$$
The only one with a root is $d_{2,Q}(t,2)$ with $t_0=2$ as a root, yielding the integral point $(2,2)$.
\end{itemize}
$C_2$: Notice that $(t_0,x_0)$ is an integral point on $C_2$ if and only if
$$x_0^3 = (x_0-1)t_0^2.$$
Noting that there are no solutions with $x_0 - 1 = 0$, we see that $x_0-1 | x_0^3.$ Since $x_0$ and $x_0-1$ share no prime factors, we must have that $x_0-1$ is $\pm 1$. So we have two possibilities for the ordered pair $(x_0,x_0-1):$
$$(x_0,x_0-1) = (2,1) \qquad \text{or} \qquad (x_0,x_0-1) = (0,-1) .$$ 
In the first case we have $8 = t_0^2$, yielding no integral (or rational) solutions. In the second, we must have $t_0 = 0$. Hence the only integral point on $C_2$ is $(t_0,x_0) = (0,0)$.

\end{proof}

Note that the $t_0 > 2$ restriction is required because $C_Q$ has the point $(2,2)$.

\begin{corollary} \label{cor:no_int_pts}
Let $t_0 > 2$ be a natural number. If the curve $C_{P}$ has no integral points of the form $(t_0,x_0)$, then the specialization map $\sigma_{t_0}$ is injective.
\end{corollary}

In order to work directly with the integral points of $C_{P}$, we will utilize the algorithm of Poulakis and Voskos \cite{Poulakis-GenusZero}. This relates finding integral points on genus zero curves to solving Pell-like equations. The algorithm requires that the number of ``valuations at infinity" (henceforth called points at infinity) of the curve is less than three; that is, there are at most two points defined over $\closure[0]{\mathbb Q}$ lying in the closure of $C_{P}$ in $\mathbb P^2$ but not on $C_{P}$ itself. Homogenizing $d_{2,P}$ then setting the new variable to zero, we obtain the equation
\begin{equation} \label{homogenized}
x^4 -4tx^3 +2t^2x^2 + 4t^3x +t^4  = 0.
\end{equation}
Setting $t=1$, we have
\begin{align*}
x^4-4x^3+2x^2-4x+1 &= 0 \\
(x^2-2x-1)^2 &= 0.
\end{align*}
On the other hand, setting $x = 1$ we similarly obtain
$$(t^2+2t-1)^2 = 0.$$
So if $\sigma$ is a root of $x^2-2x-1$ and $\tau$ is a root of $t^2+2t-1$, the points at infinity are
$$(1: \sigma : 0), (1 : \bar\sigma : 0), ( \tau: 1 : 0), (\bar\tau : 1 : 0).$$
However, notice that $1/\tau$ is a root of $x^2-2x-1$: indeed,
$$\left(\frac{1}{\tau}\right)^2 -2\frac{1}{\tau}-1 = \frac{1-2\tau-\tau^2}{\tau^2} = -\frac{\tau^2+2\tau-1}{\tau^2} = 0.$$
Hence of the four points listed above only two are distinct. Thus $C_{P}$ has two points at infinity. Poulakis and Voskos now proceed as follows.
\begin{enumerate}
	\item We first need to determine the singularities of the projective closure of $C_{P}$. Sage quickly yields $(0:0:1)$ as the only singular point.
	\item Using Sage, we obtain the rational parameterization
	$$(a:b) \to (8ab^3 + 4b^4 : 8a^2b^2 + 4ab^3 : a^4 - 4a^3b + 2a^2b^2 + 4ab^3 + b^4).$$
	Notice that the third component comes from Equation \eqref{homogenized}; in particular,
	$$a^4 - 4a^3b + 2a^2b^2 + 4ab^3 + b^4 = (a^2-2ab-b^2)^2.$$
	\item Set $u = 2a-2b$ and $v = b$. Then $a = u/2 + v$ and $b = v$. After this change of variables, our birational map becomes
	$$(u:v) \to \left(4uv^3 + 12v^4 : 2u^2v^2 + 10uv^3 + 12v^4 : \frac{1}{16}(u^2-8v^2)^2\right).$$
	Equivalently, we have
	$$(u:v) \to (16(4uv^3 + 12v^4) : 16(2u^2v^2 + 10uv^3 + 12v^4) : (u^2-8v^2)^2).$$
	Set $p(u,v) = 16(4uv^3 + 12v^4)$ and $q(u,v) = 16(2u^2v^2 + 10uv^3 + 12v^4).$
	\item The resultant $R_1$ of $p(u,1)$ and $u^2-8$ is $2^{12}$, and the resultant $R_2$ of $q(u,1)$ and $u^2-8$ is $-2^{12}$. Thus we set $D = \gcd(R_1,R_2) = 2^{12}$.
	\item Every integral point $(t_0,x_0)$ on $C_{P}$ is then obtained in the following way. Let $(u_0,v_0) \in \mathbb Z^2$ be a solution to an equation of the form $u^2-8v^2 = k$ for some $k|D$ with $u_0 \geq 0$ and $\gcd(u_0,v_0) = 1$. Then we have
	\begin{align} \label{solution_form}
	t_0 = \frac{p(u_0,v_0)}{(u_0^2-8v_0^2)^2}, && x_0 = \frac{q(u_0,v_0)}{(u_0^2-8v_0^2)^2}.
	\end{align}
\end{enumerate}	
So the specialization map $\sigma_{t_0}$ is injective for any $t_0$ which cannot be written in the form as given in \eqref{solution_form}. We will now make this even more explicit by solving the Pell-like equations given above. Many of the equations $u^2-8v^2 = k$ have no solutions of the required form, so we identify those first.

\begin{lemma} \label{lemma:Pell}
Let $l \in \mathbb{Z}$ with $l \geq 4$ or $l = 2$ and let $m \in \mathbb{Z}$ with $m \geq 4$. The equations $u^2-8v^2 = 2^l$ and $u^2-8v^2 = -2^m$ have no solutions of the form $(u_0,v_0)$ with $(u_0,v_0) \in \mathbb{Z}^2$ and gcd$(u_0,v_0) = 1$. In addition, the three equations
\begin{align*}
u^2 -8v^2 &= -2 \\
u^2 - 8v^2 &= -1, \text{ and} \\
u^2 - 8v^2 &= 2
\end{align*}
have no integer solutions at all.
\end{lemma}
\begin{proof}
Let $(u_0,v_0) \in \mathbb{Z}^2$. Suppose $u_0^2-8v_0^2 = 2^l$ with $l \geq 4$. Then $8 | u_0^2$, so necessarily $4|u_0$. Write $u_0 = 4k$ for some $k \in \mathbb{Z}$. We then have
$$2k^2-v_0^2 = 2^{l-3},$$
where $l-3 \geq 1$. Hence $2|v_0^2$, so $2|v_0$ and thus $2 | $gcd$(u_0,v_0)$, so no solutions of the required form exist.  Similarly, if $u_0^2-8v_0^2=-2^m$ for some $m \geq 4$ we also find that $2 | $gcd$(u_0,v_0)$. If $l = 2$, writing $u_0 = 2k$ we have that
$$k^2-2v_0^2 = 1.$$
Thus $2v_0^2 = (k-1)(k+1)$, so $2|k+1$ or $2|k-1$. Thus $k$ is odd, so $k-1$ and $k+1$ are both even. Hence $2|v_0^2$, so $2|v_0$ again. Thus the equations $u^2-8v^2 = 2^l$ and $u^2-8v^2=-2^m$ have no solutions $(u_0,v_0)$ with gcd$(u_0,v_0)=1$. For the remaining three equations, reducing mod 4 tells us that $u_0^2$ is congruent to either 2 or 3 mod 4, which is impossible. Thus these three equations have no integer solutions at all.
\end{proof}

Next, we show that, for the remaining equations, requiring that $\gcd(u_0,v_0) = 1$ is an extraneous condition.

\begin{lemma} 
All integer solutions $(u_0,v_0)$ to the equations $u^2-8v^2 = k$ for ${k \in \{-8,-4,1,8\}}$ have \break gcd$(u_0,v_0) = 1$.
\end{lemma}
\begin{proof}
Notice that gcd$(u_0,v_0)^2 | k$, so gcd$(u_0,v_0)$ is either 1 or 2. If gcd$(u_0,v_0)$ = 2, then setting $u_0 = 2m$ and $v_0 = 2l$ we find that $m^2-8l^2 = k/4$. If $k = 1$ then $k/4$ isn't an integer. For $k = -4$ we have ${k/4 \equiv 3 \text{ mod } 4}$ and for $k = \pm 8$ we have $k/4 \equiv 2$ mod 4. But $m^2-8l^2 \equiv$ 0 or 1 mod 4, so no $k$ allows the equality to hold mod 4. So we can't have gcd$(u_0,v_0)$ = 2, and thus gcd$(u_0,v_0) = 1$.
\end{proof}

Combining what we have shown in the previous two lemmas with step 5 from the Poulakis and Voskos algorithm, we see that the $t$-coordinates of integral points of $C_{P}$ have the form
$$t_0 = 16\frac{4u_0v_0^3+12v_0^4}{(u_0^2-8v_0^2)^2} = 64\frac{u_0v_0^3+3v_0^4}{(u_0^2-8v_0^2)^2}$$
where $(u_0,v_0)$ is an integral solution of any of the equations $u^2-8v^2 = k$ where $k \in \{-8, -4, 1, 8 \}$ and $u_0 \geq 0$. Before going any further, we use this formula for $t_0$ to extract a simple subset of $\mathbb N$ of density 1/4 for which the specialization map is injective.

\begin{theorem} \label{thm:1_mod_4} Let $t_0 \in \mathbb N$ with $t_0 > 1$ and suppose $t_0 \equiv 1$ mod $4$. Then the specialization map for $E$ at $t_0$ is injective.
\end{theorem}
\begin{proof}
Suppose that $(t_0,x_0)$ is an integral point on $C_{P}$ so that we have
$$t_0 = 64\frac{u_0v_0^3+3v_0^4}{(u_0^2-8v_0^2)^2}$$
where $u_0,v_0$ satisfies $u_0^2-8v_0^2 = m$ for some $m \in \{-8,-4,1,8\}$. Note that if $m = 1$ or $m = -4$, then $t_0$ is even. If $t_0$ is odd, we must have $m = \pm 8$, so that $t_0 = u_0v_0^3+3v_0^4$ for some $u_0,v_0$ satisfying $u_0^2-8v_0^2 = \pm 8$. Hence $8 | u_0^2$, so that $4|u_0$, and since we require that gcd$(u_0,v_0) = 1$ we have that $v_0$ is odd. Thus
$$t_0 \equiv u_0v_0^3+3v_0^4 \equiv 3 \text{ mod } 4.$$

Now assume that we have $t_0 \in \mathbb N$ with $t_0 > 1$ and $t_0 \equiv 1$ mod 4. Then we've just shown that $(t_0,x_0)$ is not an integral point on $C_P$ for any $x_0 \in \mathbb{Z}$, so by Corollary \ref{cor:no_int_pts} the specialization map at $t_0$ is injective.
\end{proof}

Using some elementary algebraic number theory, we now solve the remaining four equations.
\begin{enumerate}
\item $u^2-8v^2 = 1$: The integer solutions to this equation correspond to units of $\mathbb{Z}[\sqrt{2}]$ of the form $a+2b\sqrt{2}$ with $a,b \in \mathbb{Z}$. Recall that
$$\mathbb{Z}[\sqrt{2}]^\times = \{ \pm(1+\sqrt{2})^n \mid n \in \mathbb Z \}.$$ 
If we write $(1+\sqrt{2})^n = c+d\sqrt{2}$, note that $2|d$ if and only if $2|n$. Hence the integer solutions of $u^2-8v^2 = 1$ correspond to $\pm (1+\sqrt{2})^{2m}$ for $m \in \mathbb{Z}$. The solutions with $u \geq 0$ correspond to choosing $+$.

\item $u^2-8v^2 = -4$: Suppose $(u_0,v_0)$ is an integral solution. Noting that $u_0$ is even, set
$$x = \frac{u_0+4v_0}{2}, \hspace{.5in} y = \frac{u_0+2v_0}{2}.$$
Then $x,y$ are integers such that
$$x^2-2y^2=1;$$
that is, $x+y\sqrt{2}$ has $\mathbb Z[\sqrt{2}]$-norm 1. Note that every unit of $\mathbb Z[\sqrt{2}]$ that is an even power of $1+\sqrt{2}$ must have norm 1 because it's either a square or minus a square (and $-1$ has norm 1). Additionally, every unit of $\mathbb Z[\sqrt{2}]$ that is an odd power of $1+\sqrt{2}$ must have norm $-1$ because it's $1+\sqrt{2}$ times (plus or minus) a square, and $1+\sqrt{2}$ has norm $-1$. Hence
$$x+y\sqrt{2} = \pm(1+\sqrt{2})^{2n}$$
for some $n \in \mathbb{Z}.$

 Thus
$$u_0+4v_0+(u_0+2v_0)\sqrt{2} = \pm2(1+\sqrt{2})^{2n}.$$
Multiplying both sides by $-(1-\sqrt{2})$ gives
$$u_0+2v_0\sqrt{2} = \pm 2(1+\sqrt{2})^{2n-1}.$$
Thus the solution set of $u^2-8v^2 = -4$ corresponds to the set 
$${\{\pm 2(1+\sqrt{2})^{2n+1} \mid n \in \mathbb{Z} \} \subset \mathbb{Z}[\sqrt{2}]}.$$
 As before, the solutions with $u \geq 0$ correspond to choosing $+$.

\item $u^2-8v^2 = 8$: If $(u_0,v_0)$ is an integral solution, notice that $4|u_0$. Writing $u_0 = 4m$, we see that
$$v_0^2-2m^2 = -1.$$
Hence
$$v_0 + \frac{u_0}{4}\sqrt{2} = \pm(1+\sqrt{2})^{2n+1}.$$
Multiplying both sides by $2\sqrt{2}$, we have
$$u_0+2v_0\sqrt{2} = \pm2\sqrt{2}(1+\sqrt{2})^{2n+1}.$$
So the solution set of $u^2-8v^2 = 8$ corresponds to the set $$
\{\pm 2\sqrt{2}(1+\sqrt{2})^{2n+1} \mid n \in \mathbb{Z} \} \subset \mathbb{Z}[\sqrt{2}].$$ 
The solutions with $u \geq 0$ correspond to choosing $+$.

\item $u^2-8v^2 = -8$: As with the $k=8$ case, writing $u_0 = 4m$ we have
$$v_0^2-2m^2 = 1.$$
Using a similar argument, we find that the solution set of $u^2-8v^2 = -8$ corresponds to the set 
$$\{\pm 2\sqrt{2}(1+\sqrt{2})^{2n} \mid n \in \mathbb{Z} \} \subset \mathbb{Z}[\sqrt{2}].$$
 For $n \geq 0$ the solutions with $u \geq 0$ correspond to choosing $+$, and for $n < 0$ the solutions with $u \geq 0$ correspond to choosing $-$; notice that the choice is sgn$(n)$.
\end{enumerate}
We summarize the above discussion with the following formula that gives the $t$-coordinates of the integral points on $C_{P}$.
\begin{proposition} \label{proposition:explicit_formulas}
If $t_0$ is the $t$-coordinate of an integral point on $C_{P}$, then $t_0$ is given by one of the following four formulas.
\begin{enumerate}
\item $t_0 = 64(u_{1,n}v_{1,n}^3+3v_{1,n}^4)$ where
\begin{align*}u_{1,n} &= \frac{(1+\sqrt{2})^{2n}+(1-\sqrt{2})^{2n}}{2}, \\ v_{1,n} &= \frac{(1+\sqrt{2})^{2n}-(1-\sqrt{2})^{2n}}{4\sqrt{2}}
\end{align*}
for some $n \in \mathbb{Z}$.
\item $t_0 = 4(u_{2,n}v_{2,n}^3+3v_{2,n}^4)$ where
\begin{align*} u_{2,n} &= (1+\sqrt{2})^{2n+1}+(1-\sqrt{2})^{2n+1},\\
v_{2,n} &= \frac{(1+\sqrt{2})^{2n+1}-(1-\sqrt{2})^{2n+1}}{2\sqrt{2}}
\end{align*}
for some $n \in \mathbb{Z}$.
\item $t_0 =u_{3,n}v_{3,n}^3+3v_{3,n}^4$ where
\begin{align*}u_{3,n} &= \sqrt{2}\left((1+\sqrt{2})^{2n+1}-(1-\sqrt{2})^{2n+1}\right), \\
v_{3,n} &= \frac{(1+\sqrt{2})^{2n+1}+(1-\sqrt{2})^{2n+1}}{2}\end{align*}
for some $n \in \mathbb{Z}$.
\item $t_0 =u_{4,n}v_{4,n}^3+3v_{4,n}^4$ where
\begin{align*} u_{4,n} &= \sqrt{2}\left((1+\sqrt{2})^{2n}-(1-\sqrt{2})^{2n}\right), \\ 
v_{4,n} &= \frac{(1+\sqrt{2})^{2n}+(1-\sqrt{2})^{2n}}{2}\end{align*}
for some $n \in \mathbb{Z}$.
\end{enumerate}
\end{proposition}
\begin{proof}
Let $k=1$, so that for a solution $(u_0,v_0)$ of $u^2-8v^2=1$ we have that $t_0 = 64(u_0v_0^3+3v_0^4).$ Let $(1+\sqrt{2})^{2n} = u_0+2v_0\sqrt{2}$, so that $(1-\sqrt{2})^{2n} = u_0-2v_0\sqrt{2}.$ Adding and subtracting the equations gives
\begin{align*}
2u_0 &= (1+\sqrt{2})^{2n}+(1-\sqrt{2})^{2n} \\
4v_0\sqrt{2} &= (1+\sqrt{2})^{2n}-(1-\sqrt{2})^{2n}.
\end{align*}
Solving for the left hand sides gives the first formula, and the other 3 formulas are obtained in the exact same way. Finally, note that for formula 4, we do not need to include a sgn$(n)$ factor as discussed when solving the corresponding equation above since it cancels out in the expression for $t_0$.
\end{proof}
To summarize, we have shown the following. 
\begin{theorem} \label{theorem:density1}
Let $T$ be the set of integers $t_0 > 2$ which fail to satisfy the conditions of Proposition \ref{proposition:explicit_formulas}. Then $T \subset \mathbb N$ is a subset of density 1 and for each $t_0 \in T$ the specialization map $\sigma_{t_0}$ is injective.
\end{theorem}
\begin{proof}
Let $H$ be the Hilbert subset of $\mathbb Q$ corresponding to $d_{2,P}$. Proposition \ref{proposition:explicit_formulas} shows that $H \cap \mathbb N \subset T$. Now use Proposition \ref{prop:density1} and Corollary \ref{cor:no_int_pts}.
\end{proof}
\begin{corollary} \label{cor:ex1_rank}
Let $T$ be as in Theorem \ref{theorem:density1}. For each $t_0 \in T$, the Mordell-Weil group of the elliptic curve
$$E_{t_0} : y^2 = x^3 -t_0^2x+t_0^2$$
has a torsion-free subgroup of rank 2 generated by $(t_0,t_0)$ and $(0,t_0)$. In particular, 
$$\text{rank}(E_{t_0}(\mathbb Q)) \geq 2.$$
\end{corollary}

For our second example, we consider
$$E: y^2 = x^3 - (t^2+27)x + 10t^2+48t+90.$$
This example comes from Shioda's list of rational elliptic surfaces with specified Mordell-Weil rank \cite{Shioda-Examples}. As indicated there, $E(\mathbb Q(t))$ has rank 4 with generators
\begin{align*}
(t+3,4t+6), && (9,t+24), && (1,3t+8), && (6,2t+12).
\end{align*}
Instead of considering specialization of the entire Mordell-Weil group, we will focus on the subgroup $M$ generated by the two points $P = (t+3,4t+6)$ and $Q = (9,t+24)$ in order to show the utility of Algorithm \ref{algorithm} for proper subgroups.

As in the previous example, we obtain the four relevant polynomials.
\begin{align*}
d_{2,P}(t,x) &= x^{4} - 4 x^{3} t + 2 x^{2} t^{2} + 4 x t^{3} + t^{4} - 12 x^{3} - 68 x t^{2} - 40 t^{3} \\
 &+ 54 x^{2} - 276 x t - 258 t^{2} - 396 x - 936 t - 351 \\
d_{2,Q}(t,x) &= x^{4} + 2 x^{2} t^{2} + t^{4} - 36 x^{3} - 44 x t^{2} + 54 x^{2} - 384 x t - 306 t^{2} \\
&+ 252 x - 1728 t - 2511 \\
d_{2,P+Q}(t,x) &= x^{4} + 4 x^{3} t + 2 x^{2} t^{2} - 4 x t^{3} + t^{4} + 12 x^{3} - 92 x t^{2} + 40 t^{3} \\
& + 54 x^{2} - 492 x t + 366 t^{2} - 1044 x + 936 t + 1809 \\
g(t,x) &= x^3 - (t^2+27)x + 10t^2+48t+90
\end{align*}
Notice that the curves $C_P, C_Q$, $C_{P+Q}$ and $C_2$ have rational points 
$$(t,x) = (9,-6),(9,36), (9,6), \text{ and } (30,15),$$
respectively. Using Sage, the curves also have genus 1, so they are elliptic curves defined over $\mathbb Q$ (despite the fact that these curves are defined by quartic polynomials, there is still an embedding of their normalizations into $\mathbb P^2$ as a cubic where we move our selected rational point to infinity), and thus the methods used for the previous examples will not work. However, using Magma and Sage, we find that the curves $C_P, C_Q$ and $C_{P+Q}$ have Mordell-Weil rank zero (over $\mathbb Q$) and have the following finite lists of rational points.
\begin{align*}
	C_P(\mathbb Q) & = \{ (-11,6), (-12,9), (9,-6), (44,1) \} \\
	C_Q(\mathbb Q) & = \{ (-5,8), (-3,0), (9,36), (-1,-4) \} \\
	C_{P+Q}(\mathbb Q) & = \{(-19,-6), (-26,1), (9,6), (6,9) \} 
\end{align*}
Hence we obtain the following.
\begin{theorem} \label{thm:genus_1_ex}
Let $t_0 \in \mathbb Q$ be a rational number such that 
$$t_0 \notin \{-26,-19,-12,-11,-5,-3,-1,6,9,44\}$$
and the polynomial $g(t_0,x) = x^3 - (t_0^2+27)x + 10t_0^2+48t_0+90$ has no rational roots. Then the specialization map $\sigma_{t_0}|_M$ is injective.
\end{theorem}

\section{Acknowledgements}
I would like to thank Edray Goins and Donu Arapura for their insightful discussions regarding the contents of this paper. I also would like to thank Kenji Matsuki for contributing significant revisions and clarifications for an early draft.

This research did not receive any specific grant from funding agencies in the public, commercial, or not-for-profit sectors.

\section{Appendix}

In this appendix, we provide a Sage function \verb+HIT_spec_check+ that utilizes a variant of Algorithm \ref{algorithm} to find places where the algorithm doesn't provide useful information. The code can be directly copy-pasted into a Sage worksheet to begin computing. The code does not compute the Hilbert set $S_M$ as in Algorithm \ref{algorithm}, but instead iterates through a set \verb+test_interval+ of rational numbers and checks if the sufficient conditions outlined in \S \ref{sec:neron} are all satisfied (for $n = 2$). More specifically, the inputs are the constants $A$ and $B$ of short Weierstrass equation of an elliptic curve $E/\mathbb Q(t)$, the \verb+test_interval+, free generators of the subgroup $M$ (as in Theorem \ref{thm:neron-spec}) one would like to specialize, and the points in the torsion subgroup of $E(\mathbb Q(t))$.

 \begin{verbatim}
def list_powerset(lst): #A powerset subroutine
    result = [[]]
    for x in lst:
        result.extend([subset + [x] for subset in result])
    return result
def HIT_spec_check(A,B,test_interval,free_gens,torsion_points):
    #check for injectivity of the specialization map using Algorithm 12 for n = 2.
    #A, B Weierstrass coefficients of elliptic curve E over Q[t]
    #test_interval = list of numbers to specialize at
    #free_gens = list of generators of free part of Mordell-Weil group of E
    #torsion_points = list of 2-torsion points of E.
    #-These need to have 3 coordinates.
    print("The HIT criterion fails at:")
    for t_0 in test_interval:
        if 4*A.numerator().subs(t=t_0)^3+27*B.numerator().subs(t=t_0)^2 == 0:
            print(str(t_0) + " (not elliptic)")
            continue
        #specialized curve
        E = EllipticCurve(QQ,[A.numerator().subs(t=t_0),\
        B.numerator().subs(t=t_0)])
        #the test fails if the specialized curve gains additional 2-torsion.
        #(note: this just checks for any additional torsion currently)
        if len(E.torsion_points()) > (len(torsion_points)):
            print(str(t_0) + " (gained torsion)")
            continue
        Mmod2M_elements = []
        free_gens_spec = []
        torsion_points_spec = []
        points_to_divide = []
        for P in free_gens: #specialize the free generators
            free_gens_spec.append(E(P[0].numerator().subs(t=t_0),\
            P[1].numerator().subs(t=t_0)))
        #specialize the torsion
        for P in torsion_points:
        	#if P is the identity, we need to append the identity
            if P[2] == 0:
                torsion_points_spec.append(E(0))
            #otherwise, we specialize it
            else:
                torsion_points_spec.append(E(P[0].numerator()\
                .subs(t=t_0),P[1].numerator().subs(t=t_0)))
        #each nonempty subset corresponds to one point in Mmod2M
        powerset_free_gens = list_powerset(free_gens_spec) 
        for S in powerset_free_gens:
            #skip the empty subset
            if S == []:
                continue
            #otherwise, add up everything in the subset &
            #put it in Mmod2M_elements
            else:
                sum = E(0)
                for P in S:
                    sum = sum + P
                Mmod2M_elements.append(sum)
        #if there's no torsion, we just need to divide things in Mmod2M by 2
        if len(torsion_points) == 1:
            points_to_divide = Mmod2M_elements
        else:
            #we have to check the torsion points for divisors
            points_to_divide = torsion_points_spec 
            #this for loop adds all the other points that we have to check
            for P in torsion_points_spec: 
                for Q in Mmod2M_elements:
                    points_to_divide.append(P+Q)
        #check if any relevant points are divisible by 2
        for P in points_to_divide:
            if P.division_points(2) != []:
                print(str(t_0) + " (" + str(P) + " is divisible by 2)")
\end{verbatim}

As an example, we can obtain information about the example $E: y^2 = x^3-t^2x+t^2$. After defining the functions above, run the code below.

\begin{verbatim}
R.<t> = QQ[]
K.<t> = FunctionField(QQ)
A = -t^2
B = t^2
E = EllipticCurve([A,B])
P = E(t,t)
Q = E(0,t)
interval = QQ.range_by_height(50)
free_gens = [P,Q]
torsion_points = [(0,1,0)]
HIT_spec_check(A,B,interval,free_gens,torsion_points)
\end{verbatim}

The output generated is shown below.

\begin{verbatim}
The HIT criterion fails at:
0 (not elliptic)
1 ((-1 : -1 : 1) is divisible by 2)
-1 ((-1 : -1 : 1) is divisible by 2)
2 ((0 : 2 : 1) is divisible by 2)
-2 ((0 : -2 : 1) is divisible by 2)
3 (gained torsion)
-3 (gained torsion)
4 ((4 : 4 : 1) is divisible by 2)
-4 ((4 : 4 : 1) is divisible by 2)
1/6 (gained torsion)
-1/6 (gained torsion)
7 ((7 : 7 : 1) is divisible by 2)
-7 ((7 : 7 : 1) is divisible by 2)
8/3 (gained torsion)
-8/3 (gained torsion)
8/15 (gained torsion)
-8/15 (gained torsion)
20 ((20 : 20 : 1) is divisible by 2)
-20 ((20 : 20 : 1) is divisible by 2)
1/24 (gained torsion)
-1/24 (gained torsion)
6/25 ((0 : 6/25 : 1) is divisible by 2)
-6/25 ((0 : -6/25 : 1) is divisible by 2)
16/25 ((0 : 16/25 : 1) is divisible by 2)
-16/25 ((0 : -16/25 : 1) is divisible by 2)
27 ((27 : 27 : 1) is divisible by 2)
-27 ((27 : 27 : 1) is divisible by 2)
27/8 (gained torsion)
-27/8 (gained torsion)
27/10 (gained torsion)
-27/10 (gained torsion)
27/28 (gained torsion)
-27/28 (gained torsion)
5/49 ((-5/49 : -5/49 : 1) is divisible by 2)
-5/49 ((-5/49 : -5/49 : 1) is divisible by 2)
11/49 ((11/49 : 11/49 : 1) is divisible by 2)
-11/49 ((11/49 : 11/49 : 1) is divisible by 2)
12/49 ((-12/49 : -12/49 : 1) is divisible by 2)
-12/49 ((-12/49 : -12/49 : 1) is divisible by 2)
36/49 ((36/49 : 36/49 : 1) is divisible by 2)
-36/49 ((36/49 : 36/49 : 1) is divisible by 2)
\end{verbatim}
The output shows every rational number of height less than 50 at which Algorithm \ref{algorithm} fails to decide whether or not specialization is injective. Note, for instance, the appearance of 27 above; this is how the author discovered Example \ref{ex:linear_algebra}.

Since the code relies on finding rational points on curves defined by the 2-division polynomials of $E$, if these curves have finitely many such points (such as when they have genus larger than 1) this code can be used to conjecture the value of the constant appearing in Silverman's Specialization Theorem; that is, the ``cutoff" $c$ so that all rational numbers of height larger than $c$ have injective specialization maps.

Note that the code above could be immediately improved by readers more proficient than the author in Python and Sage. The author hopes that the code nonetheless will assist in the search for the failure of specialization maps to be injective.

\bibliography{refs}
\bibliographystyle{amsalpha}
\end{document}